\documentclass{amsart}
\usepackage{amssymb}
\newtheorem{thm}{Theorem}
\newtheorem{lem}{Lemma}

\begin{document}

\bibliographystyle{plain}

\title[Rado\v s Baki\'c]{On inequality $|z^n - 1| \geq |z-1|$}

\author[]{Rado\v s Baki\'c}

\address{U\v citeljski Fakultet, Beograd, Serbia}
\email{\rm bakicr@gmail.com}

\date{}

\begin{abstract}
We prove that $|z^n - 1| \geq |z-1|$ for all complex $z$ satisfying $|z - 1/2| \leq 1/2$ and all real $n \geq 1$.
\end{abstract}

\maketitle

\footnotetext[1]{Mathematics Subject Classification 2010 Primary 30A10.  Key words
and Phrases: Complex powers, inequalities.}

\section{Introduction}

R. Spira proved that $|w^{n+1} - 1| \geq |w^n| |w-1|$ for all Gaussian integers $w \in \mathbb C$ such that $\Re w \geq 1$ and all
positive integers $n$, see \cite{RS1}. He posed a question if this is true for all complex $w$ such that $\Re w \geq 1$, $n$ is again
a positive integer, see \cite{RS2}. The answer is affirmative, see \cite{DM}, page 140. If we set $w = 1/z$, then the inequality
is transformed into the following form: $|z^n - 1| \geq |z - 1|$ for $|z-1/2| \leq 1/2$ and $n \geq 1$ is an integer.

In this note we prove that this inequality is valid for all {\em real} $n \geq 1$.

\section{The Main Result}

As noted in the Introduction, we prove the following result.

\begin{thm}
For any real $n \geq 1$ we have
\begin{equation}
|z^n - 1| \geq |z-1|, \qquad  |z-1/2| \leq 1/2.
\end{equation}
If $n > 1$ and $z \not= 0, 1$, then the inequality is strict.
\end{thm}

Our proof of the Theorem is based on the following lemma which is also of independent interest.

\begin{lem}
If $n > 3$ then
\begin{equation}\label{lemma}
\cos^n x < 1-\sin x, \qquad \frac{2\pi}{n+1} \leq x < \frac{\pi}{2}.
\end{equation}
\end{lem}

{\it Proof of Lemma.} Taking logarithms we transform our inequality into equivalent form
\begin{equation}
\frac{n}{2} \ln (1-\sin^2 x) < \ln (1 - \sin x), \qquad  \frac{2\pi}{n+1} \leq x < \frac{\pi}{2}.
\end{equation}
Next, using power series expansion of $\ln (1-t)$ we obtain an equivalent inequality
\begin{equation}
\frac{n}{2} \sum_{k=1}^\infty \frac{\sin^{2k} x}{k} > \sum_{k=1}^\infty \frac{\sin^k x}{k}, \qquad \frac{2\pi}{n+1} \leq x < \frac{\pi}{2},
\end{equation}
or, by rearranging terms,
\begin{equation}
\sum_{k=1}^\infty \left( \frac{(n-1) \sin^{2k} x}{2k} - \frac{\sin^{2k-1} x}{2k-1} \right) > 0, \qquad \frac{2\pi}{n+1} \leq x < \frac{\pi}{2}.
\end{equation}
Clearly, it suffices to prove that each term in the sum is strictly positive, and since $\sin x > 0$ for the allowed range of values of $x$
this is equivalent to the following inequality
\begin{equation}
\sin x > \frac{2k}{(n-1)(2k-1)}, \qquad k \geq 1, \quad \frac{2\pi}{n+1} \leq x < \frac{\pi}{2}, \quad n > 3.
\end{equation}
However, the maximum of the right hand side over $k$ is attained for $k=1$ and the minimum of the left hand side over $x$ is attained
for $x = 2\pi/(n+1)$, so it suffices to verify the inequality
\begin{equation}
\sin \frac{2\pi}{n+1} > \frac{2}{n-1}, \qquad n > 3.
\end{equation}
Since $\sin x$ is strictly concave for $0 \leq x \leq \pi/2$ we have $\sin x > 2x/\pi$ for $0 < x < \pi/2$, setting $x = 2\pi/(n+1)$ this
gives
\begin{equation}
\sin \frac{2\pi}{n+1} > \frac{2}{\pi} \frac{2\pi}{n+1} > \frac{2}{n-1},
\end{equation}
the last inequality relies on the assumption $n > 3$. This proves our Lemma. $\Box$

{\it Proof of Theorem.} Since
\begin{equation}
f(z) = \frac{z-1}{z^n-1}
\end{equation}
is analytic in a neighborhood of $K = \{ z \in \mathbb C : |z-1/2| \leq 1/2 \}$ it suffices, by the Maximum Modulus Principle, to prove our inequality for $z \in \partial K = C$.  Let $z = re^{i\phi} \in C$. Since the inequality is obvious for $z = 0$ and for $z = 1$ we
can assume $0<r<1$. Since both sides of inequality are invariant under complex conjugation we can also assume $0 \leq \phi \leq \pi/2$.

If $n\phi \leq 2\pi - \phi$, then the inequality holds for elementary geometric reasons. Indeed, in that case the point $z^n$ lies
on the circle $\{ w \in \mathbb C : |w| = r^n \}$ and outside the angle $S_\phi = \{ w = \rho e^{i\theta} : \rho \geq 0, -\phi < \theta < \phi \}$. Since $r^n < r$ it is easily seen that $|z^n - 1| > |z-1|$.

Therefore, we can assume that $2\pi/(n+1) <\phi < \pi/2$. Note that this implies $n > 3$. Clearly $|z-1| = \sin \phi$, therefore we have to
prove that $z^n$ lies outside the circle $|z-1| = \sin \phi$. In fact, we are going to prove a stronger assertion: $z^n$ lies to the left
of the line $l$ given by $\Re z = 1-\sin \phi$, which is tangent to the mentioned circle. This can be expressed analytically as
$\Re z^n < 1-\sin \phi$. Since $r = \cos \phi$ this can be written as $\cos^n \phi \cos n\phi < 1-\sin \phi$. However, in view of Lemma this is immediate: $\cos^n \phi \cos n\phi \leq \cos^n \phi < 1-\sin \phi$ and the proof is finished. $\Box$

\end{document}